\documentclass[amsfonts,12pt]{article}
\def \ad{\end{document}}

  \usepackage{amssymb,amsmath}
\def \n1{\newpage}

\numberwithin{equation}{section}
\def \bqn{\begin{equation}}
\def \9{\end{equation}}
\def \3{\begin{eqnarray*}}
\def \4{\end{eqnarray*}}
\def \1{\begin{eqnarray}}
\def \2{\end{eqnarray}}

 \newcounter{corollary}

\newcounter{proposition}
\newcounter{definition}

\def \brem{\begin{remark}}
\def \erem{\end{remark}}
\def \bdf{\begin{definition}}
\def \edf{\end{definition}}
\def \bth{\begin{theorem}}
\def \eth{\end{theorem}}
\def \bprp{\begin{proposition}}
\def \eprp{\end{proposition}}
\def \bprf{\begin{proof}}
\def \eprf{\end{proof}}
\def \bex{\begin{example}}
\def \eex{\end{example}}
\def \bprf{\begin{proof}}
\def \eprf{\end{proof}}
\def \blem{\begin{lemma}}
\def \elem{\end{lemma}}
\newcounter{theorem}

\def \bl{\begin{lemma}}
\def \bcor{\begin{corollary}}
\def \ecor{\end{corollary}}
\def \el{\end{lemma}}
\def \beq*{\begin{eqnarray*}}
\def \eeq*{\end{eqnarray*}}
\def \6{\vspace*{7mm}}

\def \s1{\sqrt}

\def \bt{\begin{tabular}}
\def \et{\end{tabular}}

\def \l{\left}

\def \r{\right}
\def \hs1{\hspace*{3mm}}
\def \q2{\hspace*{9mm}}

\def \un1{\underline}

\def \vs1{\vspace*{4mm}}

\def \ba{\begin{array}}
\def \ea{\end{array}}

\newcommand{\ec}{\end{center}}
\newcommand{\bc}{\begin{center}}
\newcommand{\be}{\begin{equation}}

\newcommand{\ee}{\end{equation}}
\newcommand{\bn}{\begin{enumerate}}
\newcommand{\en}{\end{enumerate}}
\newcommand{\bi}{\begin{itemize}}
\newcommand{\ei}{\end{itemize}}
\usepackage{graphicx}

\newtheorem{theorem}{Theorem}

\newtheorem{corollary}[theorem]{Corollary}

\newtheorem{definition}[theorem]{Definition}
\newtheorem{example}{Examples}

\newtheorem{lemma}[theorem]{Lemma}

\newtheorem{proposition}[theorem]{Proposition}
\newtheorem{remark}{Remark}

\newenvironment{proof}[1][Proof]{\textbf{#1.} }{\
\rule{0.5em}{0.5em}}

\begin{document}
 \renewcommand{\thetheorem}{\thesection.\arabic{theorem}}

%\vspace*{1mm}

\title{
A NOTE ON THE JORDAN CANONICAL FORM }
\author{H. Azad\\
Department of Mathematics and Statistics\\
King Fahd University of Petroleum \& Minerals\\ Dhahran, Saudi
Arabia\\
hassanaz@kfupm.edu.sa }
\date{}
\maketitle
\baselineskip=26pt
\begin{abstract}
A proof of the Jordan canonical form, suitable for a first course
in linear algebra, is given. The proof includes the uniqueness of
the number
and sizes of the Jordan blocks. The value of the customary procedure for finding the block generators is also questioned.\\

\end{abstract}

%%
%\section{Introduction}
The Jordan form of linear transformations is an exceeding useful
result in all theoretical considerations regarding conjugacy
classes of matrices, nilpotent orbits and the Jacobson-Morozov
theorem. A classical reference for this topic is Smirnov's book
[7, p.245-254]. There is a very well known proof due to Fillipov
[3], which is also given in Strang's book [8, p.422-425]. The
American Mathematical Monthly has published at least six proofs of
the Jordan form over the years: [1, 2, 4, 5, 6, 9]. The
justification for approaching the subject yet another time can
only be the clarity and brevity of the presentation and a new
criterion for the uniqueness of the number and sizes of the Jordan
blocks. This note gives such a proof, which has the added
advantage that the most important parts can be taught in a first
course on linear algebra, as soon as basic ideas have been
introduced and the invariance of dimensions has been established.
It is thus also a contribution to the teaching of these ideas.

Although extensive work has been done in [10] regarding this
circle of ideas, the method given in this note provides a very
simple algorithm whose efficiency  is shown through  worked
examples.

In view of the algorithm given in this note, and the examples given
below, it is not clear to us why a precise determination of the
block generators is needed, although, for the sake of completeness,
we have discussed this aspect too- at the expense of increasing the
level of exposition.

It would be very desirable to compare the computational complexity
of computing the Jordan canonical form, using the algorithm given
in this paper, with other algorithms, which programmes like Maple
and Mathematica use to determine the Jordan form.

As is well known, the main technical step in establishing the
Jordan canonical form is to prove its existence and uniqueness for
nilpotent transformations. We will return to the general case
towards the end of this note.

Let $A$  be a nilpotent transformation on a finite dimensional
vector space $V$, let $v$ be a nonzero vector in $V$ and $n$ the
smallest integer such that $A^n v=0.$

\noindent {\bf Proposition 1} \textit{The vectors} $\left\{A^i
v:0\leq i<n \right\}$
\textit{are linearly independent. }\\
\bprf Take an expression
$$
\sum\limits_{i=0}^{n-1}c_i A^i v = 0, \qquad (*) $$ in which the
number of non-zero coefficients is as small as possible. If the
coefficient $c_j$  is the non-zero coefficient of largest index
$j$, then multiplying by $A^{n-j}$, we obtain an expression like
(*) of smaller length. So in (*) every $c_i$ with $i<j$ is $0$.
Therefore $c_j A^j v =0$ and therefore $A^j v =0$, with $j\leq
n-1$, which contradicts the choice of $n$. This proves the claim,
\eprf

\noindent {\bf Proposition 2}    \textit{Let} $R(A)$ \textit{be
the range space of} $A$ \textit{and} $N(A)$ \textit{be the null
space of} $A$. \textit{Let} $\left\{A(v_i ):i=1,\dots,r \right\}$
\textit{be a basis of the range space. Let} $\left\{n_j:
j=1,\dots,s \right\}$ \textit{be a basis of the null space of}
$A$. \textit{Then} $\left\{v_i: i=1,\dots,r, n_j: j=1,\dots,s
\right\}$ \textit{is a
basis of the vector space} $V$.\\
\bprf Let  $v\in V$. So $A(v)=\sum\limits_{i=1}^{r}c_i A (v_i)$.
Therefore $v-\sum\limits_{i=1}^{r}c_i v_i$   belongs to the null
space of $A$, hence it is a linear combination of the $\left\{
v_i\right\}$ and $\left\{ n_j\right\}$. To see that these vectors
are linearly independent, suppose  $\sum\limits_{i=1}^{r}c_i v_i+
\sum\limits_{j=1}^{s}d_j n_j =0$. This gives
$\sum\limits_{i=1}^{r}c A( v_i)=0$ and by linear independence of
the vectors $A(v_i)$, we get $c_i =0,\ i=1,\dots,r$. The linear
independence of $\left\{ n_j\right\}$ then shows that $d_j =0,\
j=1,\dots,s $. \eprf

\noindent {\bf Proposition 3}   $V$  \textit{is a direct sum of
cyclic subspaces.}\\
\bprf We prove this, as in the standard proofs [7,8], by induction
on dimension. The null space of $A$ is a non-zero subspace and
therefore the range space of $A$ is a proper subspace of $V$. If
this is the zero subspace, then a basis of $V$  gives the
decomposition into cyclic subspaces. So suppose that $R(A)$ is a
nonzero subspace. It is an $A$ invariant subspace. By induction on
dimensions, it is a direct sum of cyclic subspaces, with
generators  $v_i ,\ i=1,\dots,k$, and basis $A^j v_i,\ 0\leq j\leq
n_i$, and $A^{n_i +1} v_i =0$. Let $v_i=Aw_i$. So $A^{j} v_i =
AA^j w_i$  shows, using Proposition 1, that the vectors $A^{j} w_i
,\ 0\leq j\leq n_i$ are linearly independent. Also $A^{n_i +1} v_i
= A^{n_i +2} w_i =0$, so $A^{n_i +1} w_i = A^{n_i} v_i$ belong to
the null space of  $A$.

By Proposition 2, if we enlarge $A^{n_i} v_i ,\ i=1,\dots,k$, to a
basis of the null space of $A$ by adjoining independent vectors
$n_1 ,\dots,n_l$ in the null space of $A$, then $A^{j} w_i$,$0\leq
j\leq n_i$,$0\leq i\leq k$, $A^{n_i} v_i$,$i=1,\dots,k$, $n_1
,\dots,n_l$ form a basis of $V$.

Therefore, the cyclic subspaces generated by $w_i ,\ i=1,\dots,k$
and the one-dimensional subspaces generated by $n_r ,\ 1\leq r\leq
l$ give a direct sum decomposition of  $V$ into cyclic subspaces.
\eprf

From this description, it is clear that in each summand only
$A^{n_i+1} w_i = A^{n_i} v_i$ contributes to the null space of $A$
in that summand and therefore the number of summands in the above
given decomposition is the dimension of the null space of  $A$.

\noindent {\bf Corollary} Let $d_i = dim\left( N(A|R(A^i )),\
i=0,1,\dots,n \right)$, where $n$ is the smallest positive integer
so that $A^n =0$. The differences $d_0 - d_1, d_1 -
d_2,\dots,d_{n-1}-d_{n}$ give the number of Jordan blocks of sizes
$1,2,\dots,n$.\\
\bprf As shown in the proof of Proposition 3, the number of
summands in the Jordan decomposition is the dimension of the null
space of  $A$. Therefore the number of blocks of size$\geq 1$   is
$\dim(N(A))$. Applying $A$  removes all blocks, if any, of size 1,
and so the number of blocks of size$\geq 2$  is
$\dim(N(A|R(A)))=d_1$. Continuing, we get that $d_i$ is the number
of blocks of size$\geq i$, $i=1,\dots,n$. Therefore the difference
$d_{i-1}-d_{i}$ gives the number of blocks of size $i$, for
$i=1,\dots,n$. \eprf
%%%%%%%%%%%
%%
%%%%%%%%%%%%%%%%%%%%%%

 \noindent {\bf Examples} \bn

\item Let $A$  be any nilpotent upper triangular matrix whose
entries to the right of the main diagonal give a non-singular
matrix. Then the null space of $A$  is $1$ dimensional and
therefore the canonical form of $A$ consists of only one block.

In particular, the matrices
\[
\l[\ba{ccccc}
0 & 2 &  &  & \\
 & 0 & 1 &  & \\
 &  & 0 & -1 & \\
  &  &  & 0 & -2 \\
  &  &  &  & 0 \ea \r]\]
and
  \[
\l[\ba{ccccc}
0 & 1 & 2 & 3 & 4 \\
 & 0 & 7 & 6 & 5 \\
 &  & 0 & 8 & 9 \\
  &  &  & 0 & 10 \\
  &  &  &  & 0 \ea \r]\]
are conjugate matrices as they are conjugate to
\[
\l[\ba{ccccc}
0 & 1 &  &  & \\
 & 0 & 1 &  & \\
 &  & 0 & 1 & \\
  &  &  & 0 & 1 \\
  &  &  &  & 0 \ea \r].\]
%%%%%%%%%%%%%%%%%
%%
%%%%%%%%%%%%%%%%
 \item If $$  A= \left[\begin{array}{cccc}
0 & 1 & 0 & 1  \\
0 & 0 & 1 &  0 \\
0 & 0 & 0 & 1  \\
0  &0  & 0 & 0  \end{array} \right],$$  then N(A)works out to be
$1$ dimensional, so there is only $1$
Jordan block.\\
Also $$  A^{3}= \left[\begin{array}{cccc}
0 & 0 & 0 & 1  \\
0 & 0 & 0 &  0 \\
0 & 0 & 0 & 0  \\
0  &0  & 0 & 0  \end{array} \right],$$ so
$$N(A^{3})=\left[\begin{array}{c}
            x \\
            y\\
            z\\
            0\\
          \end{array}\right]$$
and, as $A^{4}=0,$ a basis of $N(A^{4})/N(A^{3})$ is
$$\nu=\left[\begin{array}{c}
            0 \\
            0\\
            0\\
            1\\
          \end{array}\right].$$
Therefore, this must be a generator of the block. \item Let $$  A=
\left[\begin{array}{cccc}
2 & 0 & 2 & 1  \\
0 & 2 & 1 &  1 \\
0 & 0 & 2 & 2  \\
0  &0  & 0 & 4  \end{array} \right]. $$  The eigenvalue $2$ is of
multiplicity $3$, so the generalized eigenspace $V_{(2)}$ is $3-$
dimensional, whose basis works out to be
$(1,0,0,0),(0,1,0,0),(0,0,1,0)$ and the
 matrix of
$A|V_{(2)}$ is therefore $$\left[\begin{array}{ccc}
  2 & 0 & 2 \\
  0 & 2 & 1 \\
  0 & 0 & 2 \\
\end{array} \right].$$
We have $$(A-2I)|V_{(2)}=\left[\begin{array}{ccc}
  0 & 0 & 2 \\
  0 & 0 & 1 \\
  0 & 0 & 0 \\
\end{array} \right].$$
Let $\tilde{A}=(A-2I)|V_{(2)}.$ This gives $d_{0}=\dim
(N(\tilde{A})=2,\,d_{1}= \dim(N(\tilde{A}|R(\tilde{A}))=1,$
$d_{2}=\dim (N(\tilde{A}|R(\tilde{A}^2))=0.$\\
Therefore $\tilde{A}$ has $d_{0}-d_{1}=1$ block of size $1$ and
$d_{1}-d_{2}=1$ block of size $2$.\\
The Jordan form of $\tilde{A}$ is therefore
$$\left[\begin{array}{ccc}
  0 & 0 & 0 \\
  0 & 0 & 1 \\
  0 & 0 & 0 \\
\end{array} \right]$$
and of $A|V_{(2)}$ is $$\left[\begin{array}{ccc}
  2 & 0 & 0 \\
  0 & 2 & 1 \\
  0 & 0 & 2 \\
\end{array} \right]$$
The eigenspace for eigenvalue $4$ is one-dimensional. Therefore,
the Jordan form of $A$ is
$$   \left[\begin{array}{cccc}
2 & 0 & 0 & 0  \\
0 & 2 & 1 &  0 \\
0 & 0 & 2 & 0  \\
0  &0  & 0 & 4  \end{array} \right]. $$ \en

%%%%%%%%%%%%%%%%%%%%
%
%     Add new stuff (Dec 6) here
%
%%%%%%%%%%%%%%%

In view of such examples, it is not clear to us why an algorithmic
procedure is needed to find the precise generators of the various
blocks, because all one needs to find the form of the Jordan blocks
is to compute the invariants $d_i$. Nevertheless, for the sake of
completeness, we outline such a procedure- at the expense of
increase in level of exposition.

\noindent \underline{Step 1:}

 Find all eigenvalues. For an eigenvalue $\lambda$,
compute the generalized eigenspace corresponding to $\lambda$.
Although, one needs to compute only all vectors annihilated by
$(A-\lambda I)^{{\dim}V}$, it is algorithmically better to compute
the vectors annihilated by $(A-\lambda I)^{n}=0$, where $n$ is the
multiplicity of the eigenvalue $\lambda$  in the characteristic
polynomial of $A$. So, by working in the generalized eigenspace for
$\lambda$, and replacing $(A-\lambda I)$ by $A$, we may assume that
$A$ is a nilpotent transformation of index $\leq n$.

From now on, we assume that $A$ is a nilpotent transformation
defined on a vector space $V$

\noindent \underline{Step 2:}

Find the number and sizes of blocks of this nilpotent transformation
according to the algorithm given below: it is a restatement of the
{\it Corollary} given on p.3. This is the most important step, which
is needed to complete the next step efficiently.

\noindent {\it  Algorithm for finding the Jordan Form}\\
For a nilpotent transformation $A$ on a finite dimensional vector
space $V$, let  $N$ be   the smallest integer such that  $A^N =0$.
Let $d_i = \dim\left( N(A|R(A^i )),\ i=0,1,\dots,N \right)$.

\noindent  The differences
$$
d_0 - d_1, d_1 - d_2,\dots,d_{N-1}-d_{N}
$$
give the number of Jordan blocks of sizes $1,2,\dots,N$.

\noindent \underline{Step 3:} {\it  Algorithm for finding the block
generators}

Call a nonzero vector $v$  is of height $n$  if $n$ is the smallest
integer so that $A^n (v)=0.$ The vector space spanned by
$v,Av,\dots,A^{n-1}v$ is $n$-dimensional. A block of size $n$ is an
$A$-invariant subspace generated by a vector of height $n$.

Let $n$  be the size of the largest block. Choose a basis of
$N(A^n)/ N(A^{n-1})$. This is a non-zero space, because there exist
blocks of size $n$. The smallest $A$-invariant subspace of the
preimages gives a direct sum of blocks, each of size  $n$. Call this
space  $W_1$.

Let $m$  be the size of the block immediately below  $n$. Consider
$N(A^m)/ N(A^{m-1})$.\\
Find a basis of $N(A^m |W_1 ) / N(A^{m-1} |W_1 ) $. \\
Let $w_1,\dots,w_r$ be the preimages of these basis elements; they
are all of height  $m$.

Extend this basis of $N(A^m |W_1 ) / N(A^{m-1} |W_1 ) $ to a basis
of $N(A^m)/ N(A^{m-1})$  by adjoining independent elements with
preimages $v_1,\dots,v_s$.

The smallest  $A$-invariant subspace spanned by $v_1,\dots,v_s$ -
call it $W_2$  has $0$ intersection with $W_1$.

Let $W_1 \oplus W_2 =W_3$.  Let $l$ be the size of the blocks, if
any, just below  $m$. Extend a basis of $N(A^l |W_3 ) / N(A^{l-1}
|W_3 ) $  to a basis of  $N(A^l)/ N(A^{l-1})$. As before, we will
get the required number of blocks of size $l$ complementary to
$W_1 \oplus W_2 =W_3$. Continuing, this will give a Jordan
decomposition.

\noindent {\it Explanation}

Step 3  is based on the following observations
\begin{enumerate}
\item If $W$ is a direct sum of blocks and the size of the
smallest block is $n$ and $0<j<n $, then the null-space of  $A^j$
in $W$ is the range space of $A^{n-j}$  in  $W$.
\item If $W$  is a direct sum of blocks of size  $n$, generated by
vectors $v_1,\dots,v_k$ - all of height  $n$, then  these vectors
are linearly independent in  $N(A^n)/ N(A^{n-1})$.
\end{enumerate}
Conversely, if vectors $w_1,\dots,w_l$  are in $N(A^n)$  and their
images in the quotient $N(A^n)/ N(A^{n-1})$ are linearly
independent, then the smallest $A$-invariant subspace generated by
$w_1,\dots,w_l$ is a direct sum of blocks of size  $n$, with
generators $w_1,\dots,w_l$.

\noindent {\bf Example}:\\
Using the above algorithm, the reader can check that if
$$
 A = \left[\begin{array}{ccccccc}
0 &1  &4  &5  &6 &7 &8  \\
 0& 0 & 1 & 6 & 7&8 & 9 \\
 0&  0&0  &0  &7 &8 & 9 \\
 0&0  & 0 & 0 &0 &10 &11  \\
 0& 0 &  0&0  &0 & 11& 12 \\
 0&0  & 0 & 0 &0 &0 & 1 \\
 0& 0 &  0& 0 &0 &0 &0
\end{array} \right]
$$
then the smallest integer $n$ so that $A^n =0$ is $6$.  There are
two blocks, of sizes $1$ and $6$ respectively, generated by

\begin{center}
$
 \left[\begin{array}{c}
0 \\
19\\
-6\\
1\\
0\\
0\\
0
\end{array} \right]
$ and $
 \left[\begin{array}{c}
0 \\
0\\
0\\
0\\
0\\
0\\
1
\end{array} \right].
$
\end{center}

%%%%%%%%%%%%%%%%%%%%
%
%     done up to this point
%
%%%%%%%%%%%%%%%

\noindent {\bf  A final remark on applications:} A main
application of the Jordan form in differential equations is in
computation of matrix exponentials. However, it is computationally
more efficient to calculate the matrix of $A$  relative to a basis
of generalized eigenvectors- not necessarily given by cyclic
vectors -and compute its exponential relative to this basis;
finally, conjugating by the change of basis matrix gives the
exponential of  $A$.

%-----------------------------------------------------------------------------%
\noindent{\bf Acknowledgments}
%-----------------------------------------------------------------------------%
The author thanks the King Fahd University of Petroleum and
Minerals for its support and excellent research facilities.

\ad

%%%%%%%%%%%%%%%%%%%%%%%%%%%%%%%%%%%%%%%%%%%%%%%%%%%%%%%%%%%%%%%%%%%%%%%%%%%%%%%%%%%%%%%%%%%%%%%%%%%%%%
%%%
%%% end of paper
%%%
%%%%%%%%%%%%%%%%%%%%%%%%%%%%%%%%%%%%%%%%%%%%%%%%%%%%%%%%%%%%%%%%%%%%%%%%%%%%%%%%%%%%%%%%%%%%%%%%%%%%%
\bibitem{} http://en.wikipedia.org/wiki/Jordan\_normal\_form